\documentclass[11pt,reqno]{amsart}
\usepackage{geometry}                % See geometry.pdf to learn the layout options. There are lots.
\usepackage{graphicx}
\usepackage{epstopdf}
\usepackage{xcolor}
\usepackage{soul}
\newtheorem{theorem}{Theorem}[section]
\newtheorem{remark}{Remark}[section]

\usepackage{amsmath,amssymb,amscd,mathrsfs,tikz-cd,mathtools,scalerel,latexsym}
\usepackage{autobreak}
\allowdisplaybreaks
\usepackage{bm}%%%%??????
\usepackage{color}
\numberwithin{equation}{section}
\DeclareMathOperator{\Htwo}{\mathcal{H}_{q^2, g/3+1}}
\DeclareMathOperator{\Atwo}{\mathcal{A}_{q^2}}
\DeclareMathOperator{\Hq}{\mathcal{H}_{q^2}}
\DeclareMathOperator{\Aq}{\mathcal{A}_{q}}
\DeclareMathOperator*{\Res}{Res}  
\newtheorem{lemma}{Lemma}[section]

\newtheorem{Coro}{Corollary}[section]
\DeclareMathOperator{\moda}{mod}
\newtheorem{definition}{Definition}[section]
\usepackage{changes}

\title[First moment of primitive quartic $L$-functions]{The first moment of central values of primitive quartic $L$-functions with fixed genus}

\author{Ziwei Hong}
\address{Mathematical School, Renmin University of China, Beijing, P.R. China}
\address{School of Mathematics and Statistics, The University of New South Wales, Sydney, Australia}
\email{hongziwei@live.com}

\begin{document}

\begin{abstract}
Let $q\equiv 3\pmod 4$ be a fixed odd prime power. We study the first moment of the central values of primitive quartic $L$-functions over $\mathbb{F}_q(T)$, averaging over characters $\chi$ of fixed genus $g$ for which $\chi^2$ is also primitive. Using a double Dirichlet series, we obtain the asymptotic formula
\begin{equation}\label{sum-abstract}
 \sum_{\substack{\chi\ \mathrm{primitive\ quartic}\\
 \chi^2\ \mathrm{primitive}\\ \operatorname{genus}(\chi)=g}}
 L_q\!\left(\tfrac12,\chi\right)
 = C_q q^{2g/3}
 + O_{q,\varepsilon}\!\left(q^{(3/5+\varepsilon)g}\right)
\end{equation}
as $g$ tends to infinity through admissible genera, where $C_q$ is an explicit constant depending only on $q$.
Thus the main term has order $q^{2g/3}$, with a power-saving error of order $q^{(3/5+\varepsilon)g}$ for any fixed $0<\varepsilon<1/15$.
Our method combines the functional equation of the individual $L$-functions with estimates for quartic Gauss sum generating functions, seeking only the region of meromorphic continuation required for the contour shift rather than a global system of functional equations for the double Dirichlet series.
\end{abstract}

\maketitle
\noindent {\bf Mathematics Subject Classification (2020)}: 11M06, 11M41, 11N37, 11L05, 11L40   

\noindent {\bf Keywords}:  central values, quartic $L$-functions, Gauss sums, function field

\section{Introduction}

  This paper investigates the mean value of Dirichlet $L$-functions $L_q(s,\chi)$ evaluated at the central point $s=\frac{1}{2}$, where $\chi$ ranges over primitive quartic Dirichlet characters of genus $g$ on $\mathbb{F}_q[T]$ with $q \equiv 3 \moda 4$.
  
 The mean value of central values of $L$-functions is a fundamental problem in analytic number theory. However, despite extensive results for quadratic and cubic $L$-functions over function fields (see \cite{andrade2012mean, florea2017improving, 2019The}) and numerous studies in the context of number fields (see \cite{berg2024vanishing, gao2021moments, dunn2025quartic, GAO2024125}), there has been relatively little attention given to the central values of quartic $L$-functions over function fields. We will use the new research paradigm provided by Gao and Zhao in \cite{GAO2024125} to fill this gap.

 In the context of number fields, Gao and Zhao \cite{gao2021moments} considered the first moment of primitive quartic $L$-functions. Using Baier and Young's method (see
 \cite{BAIER2010879}), they derived the smooth first moment:
 \begin{align}
    \sum_{(q,2)=1}\sum_{\substack{\chi \moda q\\ \chi~primitive\\  \chi^2~ remains~ primitive}}L(\frac{1}{2},\chi)\omega(\frac{q}{Q})=CQ\hat{\omega}(0)+O(Q^{9/10+\varepsilon}),
 \end{align}
 where $C$ is a positive explicit constant, $\omega$ is a smooth, compactly supported function.

 Later, Gao and Zhao \cite{GAO2024125} revisited the problem by using a novel approach, employing multiple Dirichlet series to streamline the computations. Their method significantly improved the efficiency of the original approach.

%%%%?????? 
 Our focus is the first moment in the non-Kummer function field family with fixed genus and with both $\chi$ and $\chi^2$ primitive.
 
More recently, Diaconu, Ion, Pa\c{s}ol and Popa~\cite{DIPP2026} established first-moment asymptotics for higher order Hecke characters over number fields, including the quartic case (Theorems B and E). Their method uses a system of global functional equations for multiple Dirichlet series to obtain meromorphic continuation and the bounds needed for their moment calculations. For the fixed genus non-Kummer family studied here, our approach is more direct: we seek only the region of meromorphic continuation needed for the contour shift, using the functional equation of the individual $L$-functions and bounds for quartic Gauss sum generating functions. We do not determine the full domain of convergence or construct a global system of functional equations for the double Dirichlet series. This restricted analytic setup simplifies the first moment argument for our family. The ideal indexed families in \cite{DIPP2026} differ from ours, and their function field extensions are not developed in that paper.
 %Our goal is to provide a new result over $\mathbb{F}_q[T]$.

 In the function field setting, Andrade and Keating \cite{andrade2012mean} computed the mean value of quadratic $L$-functions, subsequently improved by Florea \cite{florea2017improving}. David, Florea and Lalin \cite{2019The} computed the mean value of cubic $L$-functions.

 Both works required intricate computational techniques, relying on approximate functional equations and precise residue computations of generating functions associated with Gauss sums. Such calculations, while effective, demand substantial computational effort.

 For example, \cite{2019The} shows, for $q\equiv 2\moda 3$,
 \begin{align}\label{F}
    \sum_{\substack{\chi~ primitive~cubic\\ genus(\chi)=g}}L_q(\frac{1}{2}, \chi)=\frac{\zeta_q(\frac{3}{2})}{\zeta_q(3)}\mathcal{A}_{nK}(\frac{1}{q^2},\frac{1}{q^{3/2}})q^{g+2}+O(q^{\frac{7}{8}g+\varepsilon g}).
\end{align}
    The authors followed the standard methodology. Notably, nearly half of their exposition was the detailed derivation of an explicit formula for the residue of the generating function for cubic Gauss sums. However, this particular residue did not appear in their final asymptotic expression.

  In equation (70) of \cite{2019The}, there appears an error term of size $q^{\frac{5g}{6}}$, corresponding to the residues at $u^{3}=1$. However, this term does not occur in the final result, where the error is instead controlled by $q^{\frac{7g}{8}}$. Remark~4.5 of \cite{2019The} suggests that the $q^{\frac{5g}{6}}$ term should persist in the final expression. In contrast, \cite{hong2025meanvaluecubiclfuncitons} demonstrates that this term indeed disappears. By employing double Dirichlet series, they show that the simple poles at $u^{3}=1$ lie outside the domain of convergence, eliminating the need for explicit residue computations. Their analysis provides a theoretical justification of this fact and, in the cubic case, yields a streamlined proof precisely matching David's result.

 Motivated by \cite{GAO2024125}, \cite{2019The} and \cite{hong2025meanvaluecubiclfuncitons}, we consider the sum
 \begin{equation}\label{sum}
 \sum_{\substack{\chi\ primitive\ quartic\\ \chi^2 primitive\\ genus(\chi)=g}}L_q(\frac{1}{2}, \chi).
 \end{equation}
 Here we classify characters by their genus, which is equivalent to classifying them by the degree of their conductors. In practice one could restrict the family by imposing $\deg Q = n$ for some polynomial $Q\in\mathbb{F}_q[T]$; however, the bounds for $L$-functions such as those related to the Lindel\"of hypothesis -are governed primarily by the conductors of the underlying characters. Parameterizing the family by the modulus degree alone can therefore lead to difficulties in controlling error terms. On the other hand, this family has various conjectures by random matrix theory. See \cite{D2000RANDOM, Katz1999Zeroes}.

 For these reasons we instead parameterize the family by conductors (or equivalently by genus). The definition of the conductor of a Dirichlet character is given in Section 2.1. For the family considered here, $\chi$, $\chi^2$ and $\chi^3$ are even and have the same conductor $h$. The Riemann--Hurwitz formula for the associated quartic cyclic cover gives
 \begin{equation}
  2g=3(\deg h-2),\qquad \deg h=\frac{2g}{3}+2,
  \qquad \deg F=\frac{g}{3}+1,
 \end{equation}
where $g$ is the genus of the cover and $h=FF^\sigma$. Thus the degree of each completed nontrivial $L$-function is $2g/3$.
 
 In Section~3 of \cite{lin2024one}, Lin gave a concrete equation, showing how to feature this family. We quote it as Lemma~2.1.

 Moreover, applying Perron's formula to equation \eqref{size}, we obtain the size of this family $\ll q^{\frac{2g}{3}}$.

 The approach we will use demands only an upper bound for the quartic Gauss sum generating function $\Psi_q(f,u)$, which is more concise and direct. Our technique was first employed by Soundararajan \cite{soundararajan2010second} and later adapted by \v{C}ech \cite{Cech1} to study the ratios conjecture for real Dirichlet characters. 
 
 By employing a double Dirichlet series and applying Perron's formula, we directly rewrite the sum \eqref{sum}, bypassing the need for approximation functions or explicit residue computations. This technique provides a meromorphic continuation of the underlying double Dirichlet series to a broader domain without altering its fundamental structure, significantly reducing computational complexity. Unlike previous methods, ours does not require explicit evaluation of individual sums.

 The intended result for primitive quartic $L$-functions is stated below. In this working revision, its proof remains subject to the analytic step identified in Remark~\ref{pending-gluing}.

\begin{theorem} \label{main}
Let $q$ be an odd prime power such that $q\equiv 3\moda 4$, and let $g$ tend to infinity through admissible genera with $g/3+1\in\mathbb Z$. Then for $\varepsilon>0$,
    \begin{align}
     \sum_{\substack{\chi\ primitive\ quartic\\ \chi^2 primitive\\ genus(\chi)=g}}L_q(\frac{1}{2}, \chi)=q^{\frac{2g}{3}}(q^{2}-1)P(q^{-2})Z(q^{-2},q^{-1/2})+O(q^{(\frac{3}{5}+\varepsilon)g}),
    \end{align}
    where $P(u)$ and $Z(u,v)$ are given in Section~\ref{PZ}.
\end{theorem}

This result yields the following non-vanishing estimate.

\begin{Coro}\label{non}
    Let $q$ be an odd prime power such that $q\equiv 3\moda 4$, and let $g$ tend to infinity through admissible genera with $g/3+1\in\mathbb Z$. Then for $\varepsilon>0$,
    \begin{align*}
        \#\{\mbox{$\chi$ quartic of genus g, $\chi$ and $\chi^2$ primitive:}\ L_q(\frac{1}{2},\chi)\not=0\}\gg q^{\frac{2}{3}g-\varepsilon(\frac{g}{3}+1)}.
    \end{align*}
\end{Coro}

We restrict ourselves to the non-Kummer setting, where $q\equiv 3\moda 4$. In this scenario, the construction of primitive quartic characters is significantly more involved, requiring frequent transitions between the polynomial rings $\mathbb{F}_{q}[T]$ and $\mathbb{F}_{q^{2}}[T]$. Although our method readily extends, with minor modifications, to the Kummer case $q\equiv1\moda4$, we omit these routine extensions to maintain clarity and to emphasize the power and efficiency of our approach.

The structure of this paper is as follows: In Section 2, we introduce the notation and definitions related to primitive quartic characters, $L$-functions, and Gauss sums, along with standard results in these areas. In Section 3, we prove the analytical properties of $\Psi_q(f,u)$. Section 4 studies the generating series and its first two regions of convergence. Section 5 treats the third region, the convex hull, and the conditional contour argument, recording the analytic step that remains open.

\section{Preliminaries}
In this section, we summarize key properties of Gauss sums and $L$-functions. Since the proofs of these results are well-established and readily available in the literature, we omit them here. Readers are encouraged to focus on the quantitative relationships these properties establish, as they are crucial for the computational analysis presented later.

%\deleted{We work in the non-Kummer setting, where q is an odd prime power such that $q \equiv 3 \moda 4$. For simplicity, we exclude the Kummer case ($q \equiv 1 \moda 4$), although similar asymptotic formula for the first moment of quartic $L$-functions can be derived for that case using analogous method.}

\subsection{Primitive quartic characters and $L$-functions}

We begin by constructing primitive quartic characters of genus $g$, following the description given by Lin \cite{lin2024one}. Standard results on $L$-functions associated with quartic characters can be found in various textbooks, such as \cite{lin2024one}. For completeness, we state necessary results here.

 First, we consider the quartic character over $\mathbb{F}_{q^2}[T]$ under $q\equiv 3\moda 4$. We fix an isomorphism $\Omega$ between the fourth roots of unity in $\mathbb{C}^*$ and the quartic roots of $1$ in $\mathbb{F}^*_{q^2}$. For an irreducible monic polynomial $\pi\in\mathbb{F}_{q^2}[T]$, we define the quartic symbol $\chi_{\pi}$. Let $a\in\mathbb{F}_{q^2}[T]$. If $\pi |a$, then $\chi_{\pi}(a)=0$, and otherwise $\chi_{\pi}(a)=\alpha$, where $\alpha$ is the unique root of unity in $\mathbb{C}$ such that 
\begin{equation*}
    a^{\frac{q^{2\deg\pi}-1}{4}}\equiv\Omega(\alpha)\moda \pi.
\end{equation*}
We extend the definition by multiplicity. For any monic polynomial $F\in\mathbb{F}_{q^2}[T]$, $F=\pi_1^{e_1}\dots\pi_s^{e_s}$ with distinct $\pi_i$, we define $\chi_F=\chi_{\pi_1}^{e_1}\dots\chi_{\pi_s}^{e_s}$.

Then,  we consider the primitive characters over $\mathbb{F}_{q}[T]$ with $q\equiv 3\moda 4$. \cite{2017On} gives the following description of the family: $\chi$ is quartic and both $\chi$ and $\chi^2$ are primitive over $\mathbb{F}_q[T]$ if and only if $\chi=\chi_{F}|_{\mathbb{F}_q[T]}$ where $F\in\mathbb{F}_{q^2}[T]$ is square-free and $F$ has positive degree and no nonconstant divisor in $\mathbb{F}_q[T]$. In this case, $\chi_F$ has conductor $FF^{\sigma}$, where $\sigma$ is the generator of Gal$(\mathbb{F}_{q^2}/\mathbb{F}_q)$. For more details, see \cite{2017On}. 

% We restrict $\chi_{F}$ to $\mathbb{F}_q[T]$, 
%Specifically, let $\chi_F(\cdot)=\left(\frac{\cdot}{F}\right)_3$ denote the cubic residue symbol on $\mathbb{F}_{q^2}[T]$, where $F\in\mathbb{F}_q[T]$ is a monic, square-free polynomial. Thus, it is natural to view $\chi_F(\cdot)$ as a cubic character on $\mathbb{F}_q[T]$.  

 Let $\mathcal{A}_q$ denote the set of monic polynomials over $\mathbb{F}_q$. Let $\Htwo$ represent the set of monic, square-free polynomials in $\Atwo$ of degree $g/3+1$. We define the modulus $|a|=q^{\deg a}$ for $a\in\Aq$ and $|a|_2=q^{2\deg a}$ for $a\in\Atwo$. With this notation, we have
 \begin{lemma}[Lemma~3.6 in \cite{lin2024one}]
      Suppose $q\equiv 3\moda 4$. Then
     \begin{equation}\label{equation}
 \sum_{\substack{\chi\ primitive\ quartic\\ \chi^2 primitive\\
 genus(\chi)=g}}L_q(\frac{1}{2}, \chi)=  \sum_{\substack{F\in\Htwo \\ (F,F^\sigma)=1}}L_q(\frac{1}{2}, \chi_{F}).
 \end{equation}
 \end{lemma}

For $a\in\mathbb F_q^\times$ and a prime $\pi\in\mathbb F_{q^2}[T]$ of degree $d$, the integer $(q^{2d}-1)/4$ is divisible by $q-1$, since $4\mid q+1$. Hence $\chi_\pi(a)=1$, and every character in our family is even. We therefore have the following functional equation:
 \begin{lemma}[Functional equation for even $L$-functions] \label{fe}
     Let $F\in\mathcal{H}_{q^2}$ have positive degree and satisfy $(F,F^\sigma)=1$. Then
     $$L_q(s,\chi_{F})=\omega(\chi_F)q^{2s-1}\frac{1-q^{-s}}{1-q^{s-1}}\frac{L_q(1-s,\overline{\chi_F})}{|F|_2^{s-\frac{1}{2}}},$$
     where $\omega(\chi_F)=-q^{1-\deg F}\sum\limits_{\substack{f\in\Aq\\ \deg f=2\deg F-1}}\chi_F(f)$.
 \end{lemma}
 \begin{proof}
    The general functional equation for a primitive even character of modulus h is given by
     \begin{align}
         L_q(s,\chi)=\omega(\chi)(q^{\frac{1}{2}}q^{-s})^{\deg h-2}\frac{1-q^{-s}}{1-q^{s-1}}L_q(1-s,\Bar{\chi}).
     \end{align}
     In our case, $\chi=\chi_F$ has conductor of degree $2\deg F$, which proves the formula.
 \end{proof}
 %Here $|F|_2=q^{2\deg F}$ represents the modulus of polynomials in $\mathbb{F}_{q^2}[T]$.

 For any odd character $\chi$ on $\mathbb{F}_q[T]$ (or $\mathbb{F}_{q^2}[T]$), we define $\tau(\chi)$ as the Gauss sum of the restriction of $\chi$ to $\mathbb{F}_q^*$ (or $\mathbb{F}_{q^2}^*$), i.e.,
 \begin{align*}
     \tau(\chi)=\sum_{a\in\mathbb{F}_q^*}\chi(a)e^{2\pi i \text{Tr}_{\mathbb{F}_q/\mathbb{F}_p}(a)/p}
 \end{align*}
 (or $\tau(\chi)=\sum\limits_{a\in\mathbb{F}_{q^2}^*}\chi(a)e^{2\pi i tr_{\mathbb{F}_{q^2}/\mathbb{F}_p}(a)/p}$).
 We define the sign of the Gauss sum by 
 \begin{align*}
     \epsilon(\chi)=q^{-1/2}\tau(\chi).
 \end{align*}
 When $\chi$ is even, we set $\epsilon(\chi)=1$.
\begin{Coro}[Corollary~2.4 in \cite{2019The}]\label{relation}
    Let $\chi$ be a primitive character of modulus $h$ on $\mathbb{F}_q[T]$. Then
    \begin{align*}
        \omega(\chi)=
        \begin{cases}
           \frac{1}{\tau(\chi)}q^{-\frac{\deg h-1}{2}}G_q(\chi),  &  \mbox{if $\chi$ is odd;}\\
          q^{-\frac{\deg h}{2}}G_q(\chi),   & \mbox{if $\chi$ is even.}
        \end{cases}
    \end{align*}
\end{Coro}

 This is a standard result on functional equation for $L$-functions, and its proof can be found in \cite{2019The}. Here $G_q(\chi)$ is the Gauss sum defined by
  \begin{definition}
     For $\chi$ a primitive character of the modulus h on $\mathbb{F}_q[T]$, let
     $$G_q(\chi)=\sum_{a \moda h}\chi(a)e_q(\frac{a}{h}),$$
     where $$e_q(a/h)=e^{2\pi i\text{Tr}_{\mathbb{F}_q/\mathbb{F}_p}(a_1)/p}$$ where $a_1$ is the coefficient of $1/T$ of $\frac{a}{h}$ in $\mathbb{F}_q((1/T))$.
 \end{definition}
 
 Next, we present the upper bound of $L$-functions. The following lemma restates Lemma 2.6 and Lemma 2.7 from \cite{2019The}.
 \begin{lemma}[Lindel\"of Hypothesis]
 \label{LLH}
     Let $\chi$ be a nonprincipal primitive Dirichlet character of conductor $h$ over $\mathbb{F}_q[T]$. The standard bounds below apply in particular to the quadratic and quartic characters used here. Then, for $\Re(s)\ge\frac{1}{2}$ and all $\varepsilon>0$,
     \begin{align}
         |L_q(s,\chi)|\ll q^{\varepsilon\deg h};
     \end{align}
     for $\Re(s)\ge 1$ and for all $\varepsilon>0$
     \begin{align}
         |L_q(s,\chi)|\gg q^{-\varepsilon\deg h}.
     \end{align}
 \end{lemma}
%In our case, the primitive cubic character $\chi_F$ has conductor $F$.
%Combing functional equation, we have
% \begin{lemma}\label{bound for L}
 %   Let $\chi$ be a primitive cubic character of conductor $h$ defined over $\mathbb{F}_q[T]$. Then, for all $\varepsilon>0$,
 %   \begin{align}
  %      |L_q(s,\chi)|\ll\begin{cases} q^{\varepsilon\deg h},& \mbox{for}~ \Re(s)\ge\frac{1}{2};\\
   %         \left|\frac{q^{-2s-1}-q^{-3s-1}}{1-q^{s-1}}\right|q^{\deg h(\frac{1}{2}+\varepsilon-\Re(s))},& \mbox{otherwise.}
    %    \end{cases}
   % \end{align}
% \end{lemma}
% \begin{proof}
%Combined Lemma~\ref{fe} and Lemma~\ref{LLH}, we obtain
%\begin{align*}
%    |L_q(s,\chi)|\ll \left|q^{-2s-1}\frac{1-q^{-s}}{1-q^{s-%1}}\right||h|^{\frac{1}{2}+\varepsilon-\Re(s)}
%\end{align*}
%for $\Re(s)<\frac{1}{2}$. Here $|h|=q^{\deg h}$ stands for modules of polynomials in $\mathbb{F}_q[T]$.
% \end{proof}

 %Detailed discussion of cubic characters on $\mathbb{F}_q[T]$ can be found in \cite{2019The,2017On}.

\subsection{Gauss sums}
 
We recall some well-known properties of polynomial Gauss sums. For further details, see \cite{ZHENG2017460}.

 We define generalized Gauss sum as follows.
 \begin{definition}
 Let $\chi_f$ be character of $\mathbb{F}_q[T]$ modulo $f$ for some $f\in\mathbb{F}_q[T]$. We define generalized Gauss sum
     \begin{align}
          G_q(V, f)=\sum_{u\moda f}\chi_f(u)e_q\left(\frac{uV}{f}\right).
     \end{align}
 \end{definition}
 If $(a,f)=1$, then
 $$G_q(aV,f)=\Bar{\chi}_f(a)G_q(V,f).$$
 
\begin{lemma}[Quartic Gauss sums under restriction]\label{quartic-restriction}
Suppose that $q\equiv3\pmod4$, and put
\[
 \kappa_q=(-1)^{(q+1)/4}
 =\begin{cases}-1,&q\equiv3\pmod8,\\1,&q\equiv7\pmod8.\end{cases}
\]
Let $F\in\mathcal H_{q^2}$ have positive degree and satisfy $(F,F^\sigma)=1$, where $\sigma$ is the nontrivial automorphism of $\mathbb F_{q^2}/\mathbb F_q$. Then
\begin{equation}\label{quartic-gauss-transfer}
 G_q(\chi_F)=\chi_F(F^\sigma)G_{q^2}(1,F)
 =\kappa_q^{\deg F}G_{q^2}(1,F).
\end{equation}
For a square-free $D\in\mathcal A_q$ and $(N,D)=1$, we also have
\begin{equation}\label{gauss-D-bound}
 |G_{q^2}(N,D)|=q^{\deg D}.
\end{equation}
\end{lemma}
\begin{proof}
Write $h=FF^\sigma$. The Chinese remainder theorem parametrizes the residue classes modulo $h$ over $\mathbb F_q$ by
\[
 a=\beta F^\sigma+\beta^\sigma F\pmod h,
 \qquad \beta\pmod F.
\]
For these representatives,
\[
 \chi_F(a)=\chi_F(\beta)\chi_F(F^\sigma),\qquad
 e_q(a/h)=e_{q^2}(\beta/F).
\]
Summing over $\beta$ proves the first equality in \eqref{quartic-gauss-transfer}.
Let $n=\deg F$ and $R=\operatorname{Res}(F,F^\sigma)\ne0$. The definition of the quartic residue symbol, together with the norm formula for the resultant, gives
\[
 \Omega\bigl(\chi_F(F^\sigma)\bigr)=R^{(q^2-1)/4}.
\]
Since
\[
 R^q=\operatorname{Res}(F^\sigma,F)=(-1)^{n^2}R=(-1)^nR,
\]
we obtain
\[
 R^{(q^2-1)/4}=(R^{q-1})^{(q+1)/4}
 =(-1)^{n(q+1)/4}=\kappa_q^n.
\]
The isomorphism $\Omega$ fixes $\{1,-1\}$, proving the second equality.
Finally, a square-free $D$ remains square-free over $\mathbb F_{q^2}$, and its quartic residue character is primitive modulo $D$. The usual absolute value formula for a primitive Gauss sum, and the change of variable for $(N,D)=1$, give
$|G_{q^2}(N,D)|=|D|_2^{1/2}=q^{\deg D}$ (with $D=1$ interpreted directly).
\end{proof}

\begin{lemma}
\label{Gauss1} 
    Suppose that $q\equiv 1 (\moda 8)$.
    \begin{enumerate}
        \item If $(f_1,f_2)=1$, then
        \begin{align*}
            G_q(V,f_1f_2)=&\chi_{f_1}(f_2)^2G_q(V,f_1)G_q(V,f_2)\\
            =&G_q(Vf_2^2,f_1)G_q(V,f_2).
        \end{align*}
        \item If $V=V_1P^{\alpha}$ where $P\not| V_1$, then
        \begin{equation}
            G_q(V,P^i)=\left\{
            \begin{tabular}{ll}
                0, & if $i\le\alpha$ and $i\not\equiv 0 \moda 4$;\\
                $\phi(P^i)$, & if $i\le\alpha$ and $i\equiv 0 \moda 4$;\\  
                $-|P|^{i-1}$, & if $i=\alpha+1$ and $i\equiv 0\moda 4$;\\              $\epsilon(\chi_P^i)\omega(\chi_P^i)\chi_P^i(V_1^{-1})|P|^{i-\frac{1}{2}}$, & if $i=\alpha+1$ and $i\not\equiv 0\moda 4$;\\
                0, & if $i\ge\alpha+2$,            \end{tabular}\right.
        \end{equation}
        where $\phi$ is the Euler function for polynomials. For $4\nmid i$, the character $\chi_P^i$ in the root number is the nonprincipal primitive character modulo $P$. We use $\epsilon(\chi)=\tau(\chi)/|\tau(\chi)|$ for odd $\chi$, and $\epsilon(\chi)=1$ for even $\chi$. 
    \end{enumerate}
\end{lemma}

\begin{remark}
    If we have $q=4k+3$, then $q^2=16k^2+24k+9\equiv 1\moda 8$. 
\end{remark}

\begin{proof}
    The first statement is a standard result. We write $u \moda f_1f_2$ as $u=u_1f_1+u_2f_2$ for $u_1 \moda f_2$ and $u_2 \moda f_1$, Then
    \begin{align*}
        G_q(V,f_1f_2)=&\chi_{f_2}(f_1)\chi_{f_1}(f_2)\sum_{u_1\moda f_2}\sum_{u_2\moda f_1} \chi_{f_1}(u_2)\chi_{f_2}(u_1)e_q\left(\frac{u_1V}{f_2}\right)e_q\left(\frac{u_2V}{f_1}\right)\\
        =&\chi_{f_1}(f_2)^2G_q(V,f_1)G_q(V,f_2).
    \end{align*}
    For $i\le\alpha$, the additive character is trivial, so
    \begin{align*}
      G_q(V,P^i)&=\sum_{u\bmod P^i}\chi_{P^i}(u)
      =|P|^{i-1}\sum_{u\bmod P}\chi_P(u)^i\\
      &=\begin{cases}\phi(P^i),&4\mid i,\\0,&4\nmid i.\end{cases}
    \end{align*}
    This follows from character orthogonality on $(\mathbb F_q[T]/P)^\times$.

    Now assume that $i=\alpha+1$. Write $u \moda P^i$ as $u=PA+C$ with $A\moda P^{i-1}$ and $C\moda P$. Then
    \begin{align*}
        G_q(V,P^i)&=\sum_{A\moda P^{i-1}}\sum_{C\moda P}\chi_{P^{i}}(C)e_q\left(\frac{CV_1}{P}\right)\\ & =|P|^{i-1}\chi_{P^i}(V_1^{-1})\sum_{C\moda P}\chi_{P^i}(C)e_q\left(\frac{C}{P}\right).
    \end{align*}
    If $i\equiv 0\moda 4$, then 
    \begin{align*}
        G_q(V,P^i)=|P|^{i-1}\sum_{C\not=0}e_q(\frac{C}{P})=-|P|^{i-1}.
    \end{align*}
    If $i\not\equiv0\pmod4$, the character $\chi_P^i$ is primitive modulo $P$, and the preceding identity becomes
    \[
      G_q(V,P^i)=|P|^{i-1}\chi_P^i(V_1^{-1})G_q(\chi_P^i).
    \]
    Corollary~\ref{relation}, applied to $\chi_P^i$, gives the stated formula.
    For $i\ge \alpha+2$, see \cite{florea2017improving}.
\end{proof}

\subsection{Multivariable complex analysis lemmas}
Our approach relies on two foundational results from multivariable complex analysis. We begin by introducing the concept of a tube domain.
\begin{definition}
		An open set $T\subset\mathbb{C}^n$ is a tube if there is an open set $U\subset\mathbb{R}^n$ such that $T=\{z\in\mathbb{C}^n:\ \Re(z)\in U\}.$
\end{definition}
	
   For any set $U\subset\mathbb{R}^n$, we define $T(U)=U+i\mathbb{R}^n\subset \mathbb{C}^n$.  We quote the following Bochner's Tube Theorem \cite{Boc}.
\begin{theorem}
\label{Bochner}
		Let $U\subset\mathbb{R}^n$ be a connected open set and $f(z)$ be a function holomorphic on $T(U)$. Then $f(z)$ has a holomorphic continuation to
the convex hull of $T(U)$.
\end{theorem}

 We denote the convex hull of an open set $T\subset\mathbb{C}^n$ by $\widehat T$.  Our next result is \cite[Proposition C.5]{Cech1} on the modulus of holomorphic continuations of multivariable complex functions.
\begin{theorem}

\label{Extending inequalities}
		Assume that $T\subset \mathbb{C}^n$ is a tube domain, $g,h:T\rightarrow \mathbb{C}$ are holomorphic functions, and let $\tilde g,\tilde h$ be their
holomorphic continuations to $\widehat T$. If  $|g(z)|\leq |h(z)|$ for all $z\in T$ and $h(z)$ is nonzero in $T$, then also $|\tilde g(z)|\leq
|\tilde h(z)|$ for all $z\in \widehat T$.
\end{theorem}

\section{Averages of Gauss sums}

%\subsection{Generating function of quartic Gauss sums}
In this section, assume $q\equiv 1\moda 8$. We define the generating functions for quartic Gauss sums
\begin{align}
    \Psi_q(f,u)=\sum_{F\in\Aq}G_q(f,F)u^{\deg F}
\end{align}
and
\begin{align}
    \Tilde{\Psi}_q(f,u)=\sum_{\substack{F\in\Aq\\ (f,F)=1}}G_q(f,F)u^{\deg F}.
\end{align}

Let $\eta \in \mathbb{F}_q((1/T))^\times$ and define
\[\psi(f,\eta, u)=(1-u^4q^{4})^{-1} \sum_{\substack{F\in \Aq\\ F \sim \eta}} G_q(f, F)u^{\deg(F)},\]
where the equivalence relation is given by 
\[F\sim \eta \Leftrightarrow F/\eta \in (\mathbb{F}_q((1/T))^\times)^4.\]

There is a distinction between our definition of  $\psi(f,\eta, u)$ and the one in \cite[p. 245]{https://doi.org/10.1112/plms/s3-54.2.193}. Here we sum over monic polynomials in $\mathbb{F}_q[T]$, whereas in 
\cite{1992Theta} the sum is over all polynomials in $\mathbb{F}_q[T]$. The difference accounts for extra factors of the form $(q-1)/4$ in \cite{https://doi.org/10.1112/plms/s3-54.2.193}.
Since we are considering monic polynomials, it suffices to consider the equivalence classes 
that separate degrees, i.e. $\eta=\pi_\infty^{-i}$, where $\pi_\infty$ is the uniformizer of the prime at infinity, i.e. $T^{-1}$ in $\mathbb{F}_q((1/T))$.

Basic knowledge of algebra in $\mathbb{F}_q((1/T))$ shows that for any $i \in \mathbb{Z}$, 
\[\psi(f,\pi_\infty^{-i}, u)=(1-u^4q^{4})^{-1} \sum_{\substack{F\in \Aq\\ \deg F \equiv i \moda 4}} G_q(f, F)u^{\deg F}.\]

Thus $\psi(f,\pi_\infty^{-i}, u)$ depends only on the value of $i$ modulo $4$. 

We can express the generating series $\Psi_q(f,u)$ as
\begin{eqnarray} \label{gen-series}
 \Psi_q(f,u)= (1-u^4q^4) \sum_{i=0}^{3} \psi(f, \pi_\infty^{-i},u).\end{eqnarray}

Taking $a=1$ and $l=4$ in Corollary~5.4 of \cite{david2025nonvanishing}, we obtain, for quartic Gauss sums,
\begin{lemma}\label{upper2}
Let $u=q^{-s}$ and $\sigma=\Re(s)$. For $1\le\sigma<\frac32$, equivalently $q^{-3/2}<|u|\le q^{-1}$, on compact subsets avoiding the possible poles, we have
\begin{align}\label{psibar}
\tilde{\Psi}_q(f,u)\ll |f|^{\frac{1}{2}(\frac{3}{2}-\sigma)+\varepsilon}.
\end{align}
\end{lemma}

\section{Analytic behavior of $A_4(u,v)$}
  \subsection{The generating function}
 Now, we construct a double Dirichlet series to rewrite the sum in equation \eqref{equation}. From now on, we will always use $P$, $P_1$ or $P_2$ to denote monic irreducible polynomials in the corresponding polynomial ring. And we assume $q\equiv 3\moda 4$.
 \begin{remark}
     In the previous section we worked over the base field $\mathbb{F}_{q}$ under the assumption $q \equiv 1 \moda 8$. In this section we instead consider the base field $\mathbb{F}_{q^{2}}$ with $q \equiv 3 \moda 4$. Since $q^{2} \equiv 1 \moda 8$, the arguments from the preceding section carry over verbatim: replacing $\mathbb{F}_{q}$ by $\mathbb{F}_{q^{2}}$ yields the corresponding results.
 \end{remark}
 For $\Re(s)$ and $\Re(w)$ sufficiently large, we define
 \begin{equation}
 \label{A3}
     A_4(s, w)=\sum_{\substack{F\in\mathcal{H}_{q^2} \\ (F,F^\sigma)=1}} \frac{L_q(w, \chi_{F})}{|F|^s_2}.
 \end{equation}

The term $F=1$ contributes $\zeta_q(v)=(1-qv)^{-1}$, where $v=q^{-w}$. We retain $A_4$ for the Euler product calculations and set
\begin{equation}\label{A-positive}
 \widehat A_4(u,v)=A_4(u,v)-\frac{1}{1-qv}
 =\sum_{\substack{F\in\mathcal H_{q^2}\\\deg F\ge1,\ (F,F^\sigma)=1}}
 L_q(v,\chi_F)u^{\deg F}.
\end{equation}
In additive coordinates, $\widehat A_4(s,w)$ denotes the same function with $u=q^{-2s}$ and $v=q^{-w}$. The subtraction leaves every positive-degree coefficient in $u$, and the residue at $u=q^{-2}$, unchanged. It removes the extraneous pole at $v=q^{-1}$ from the continuation argument.
 By Perron's formula for function fields, we have
 \begin{equation}
     \sum_{n\le N}a(n)=\frac{1}{2\pi i}\oint_{|u|=r}(\sum_{n=0}^{\infty}a(n)u^n)\frac{du}{(1-u)u^{N+1}},
 \end{equation}
which allows us to rewrite the sum as 
\begin{align}
    \sum_{\substack{F\in\Htwo \\ (F,F^\sigma)=1}}L_q(\frac{1}{2}, \chi_{F})
    =& \frac{1}{2\pi i}\oint_{|u|=r}\sum_{\substack{F\in\mathcal{H}_{q^2} \\ (F,F^\sigma)=1}} L_q(\frac{1}{2}, \chi_{F}) u^{\deg F}\frac{du}{u^{g/3+2}}\notag\\
    =&\frac{1}{2\pi i}\oint_{|u|=r}A_4(u,\frac{1}{2})\frac{du}{u^{g/3+2}},
\end{align}
where $A_4(u, w)=A_4(s, w)$ upon taking $u=q^{-2s}$. Hereafter, we will use $A_4(s,w)$ and $A_4(u,w)$ interchangeably when there is no ambiguity.

Next, we consider the different expressions of $A_4(s,w)$ on different convergence regions and find out the convex hull of convergence regions.

 \subsection{The first convergence region}
 To analyze the behavior of $A_4(u,v)$, we begin by using M\"obius inversion 
 \begin{equation}     \sum_{\substack{D\in\Aq\\ D|F}}\mu(D)=\left\{ \begin{array}{ll}
         1 & \mbox{$F$ has no prime divisor in } \mathbb{F}_q[T]; \\
         0 & \mbox{otherwise}.
     \end{array}
     \right.
 \end{equation}
 This allows us to remove the divisor condition.
 Applying this, $A_4(s,w)$ can be rewritten as
  \begin{align}
  A_4(s, w)=&\sum_{\substack{F\in\Hq \\ (F,F^\sigma)=1}}\frac{L_q(w,   \chi_{F})}{|F|_2^s}\notag\\
  %=&\sum_{F\in\Hq}\sum_{\substack{D|F\\ D\in\Aq}} \mu(D)\frac{L_q(w, \chi_F)}{|F|^s}\notag\\
  =&\sum_{D\in\Aq} \mu(D)\sum_{\substack{F\in\mathcal{H}_{q^2}\\ (D,F)=1}} \frac{L_q(w, \chi_{DF})}{|DF|_2^s}\notag\\
  =&\sum_{D\in\Aq}\mu(D)\sum_{\substack{F\in\mathcal{H}_{q^2}\\ (D,F)=1}}\sum_{N\in\Aq}\frac{\chi_{DF}(N)}{|N|^{w}|DF|_2^s}\notag\\
  =&\sum_{N\in\Aq}\frac{1}{|N|^w}\sum_{\substack{D\in\Aq\\ }}\frac{\mu(D)\chi_D(N)}{|D|_2^s}\sum_{\substack{F\in\Hq\\ (D,F)=1}}\frac{\chi_F(N)}{|F|_2^s}\notag\\
  =&\sum_{N\in\Aq}\frac{1}{|N|^w}\frac{L_{q^2}(s, \chi^{(N)})}{L_{q^2}(2s, \chi^{2(N)})}\sum_{\substack{D\in\Aq\\ }}\frac{\mu(D)\chi_D(N)}{|D|_2^s}\prod_{\substack{P_2|D\\ P_2\in\Atwo}}\left(1+\frac{\chi_{P_2}(N)}{|P_2|_2^s}\right)^{-1}.
 \end{align}
 The last equation uses a trick that
 \begin{align*}
     \sum_{\substack{F\in\Hq\\ (D,F)=1}}\frac{\chi_F(N)}{|F|_2^s}&=\prod_{\substack{P_2\in\Atwo\\ P_2\nmid D}}\left(1+\frac{\chi_{P_2}(N)}{|P_2|^s_2}\right)\\ & =\prod_{\substack{P_2\in\Atwo}}\left(1+\frac{\chi_{P_2}(N)}{|P_2|^s_2}\right)\prod_{\substack{P_2\in\Atwo\\ P_2| D}}\left(1+\frac{\chi_{P_2}(N)}{|P_2|^s_2}\right)^{-1}\\
     & =\frac{L_{q^2}(s,\chi^{(N)})}{L_{q^2}(2s,\chi^{2(N)})}\prod_{\substack{P_2\in\Atwo\\ P_2| D}}\left(1+\frac{\chi_{P_2}(N)}{|P_2|^s_2}\right)^{-1}.
 \end{align*}
 Here $\chi^{(N)}(F)=\chi_F(N)$ is a multiplicative character induced by $\chi_F$. $L_{q^2}(s,\chi^{(N)})$ denotes the $L$-function over $\mathbb{F}_{q^2}(T)$. %We recall that $|F|=q^{\deg F}$ and $|F|_2=q^{2\deg F}$.

 Using the lemma below, we simplify $A_4(s,w)$ further.
 \begin{lemma}
     If $D, F\in\Aq$ and $(D, F)=1$, then we have $\chi_D(F)=1$.
 \end{lemma}
 \begin{proof}
     We suppose $D$ is irreducible in $\Aq$ and show that $\chi_D(F)=1$ for all prime $D\in\Aq$. If $2|\deg D$, then $D=\pi\pi^{\sigma}$ for some prime $\pi\in\Atwo$, where $\sigma$ is the generator of Gal$(\mathbb{F}_{q^2}/\mathbb{F}_{q})$. It's easy to see that $\chi_{\pi^{\sigma}}|_{\Aq}=\overline{\chi_{\pi}|_{\Aq}}$. Hence, we must have $\chi_D(F)=\chi_{\pi}(F)\chi_{\pi^{\sigma}}(F)=1$ for $F\in\Aq$ and $(F,D)=1$. If $d=\deg D$ is odd, then $D$ remains irreducible over $\mathbb F_{q^2}$. The residue of $F$ lies in $\mathbb F_{q^d}^{\times}$, and
\[
 \frac{q^{2d}-1}{4}=(q^d-1)\frac{q^d+1}{4}
\]
is a multiple of $q^d-1$, since $q^d\equiv3\pmod4$. Thus $\chi_D(F)=1$ in this case as well. Multiplicativity proves the assertion for arbitrary $D$.
     
 \end{proof}
 Using the lemma, we obtain
\small{
\begin{align*}
    A_4(s, w)=&\sum_{N\in\Aq}\frac{1}{|N|^{w}}\frac{L_{q^2}(s, \chi^{(N)})}{L_{q^2}(2s, \chi^{2(N)})}\sum_{\substack{D\in\Aq\\ (D,N)=1}}\frac{\mu(D)}{|D|^s_2}\prod_{\substack{P_2\in\Atwo\\ P_2|D}}\left(1+\frac{\chi_{P_2}(N)}{|P_2|^s_2}\right)^{-1}\\
    =&\sum_{N\in\Aq}\frac{1}{|N|^{w}}\frac{L_{q^2}(s, \chi^{(N)})}{L_{q^2}(2s, \chi^{2(N)})}\prod_{\substack{P_1\in\Aq\\ (P_1,N)=1}}\left(1-\frac{1}{|P_1|^s_2}\prod_{\substack{P_2|P_1\\
    P_2\in\Atwo}}\left(1+\frac{\chi_{P_2}(N)}{|P_2|_2^s}\right)^{-1}\right) \\
    =&\sum_{N\in\Aq}\frac{1}{|N|^{w}}\frac{L_{q^2}(s, \chi^{(N)})}{L_{q^2}(2s, \chi^{2(N)})}P(s,\chi^{(N)})\prod_{\substack{P_1|N\\ P_1\in\Aq}}\left(1-\frac{1}{|P_1|_2^s}\right)^{-1},
\end{align*}}
where
\begin{equation}
    P(s, \chi^{(N)})=\prod_{P_1\in\Aq}\left(1-\frac{1}{|P_1|_2^s}\prod_{\substack{P_2|P_1\\ P_2\in\Atwo}}\left(1+\frac{\chi_{P_2}(N)}{|P_2|_2^s}\right)^{-1}\right).
\end{equation}

Write $\sigma=\Re(s)>\frac12$. Uniformly in $N$, the local factors give 
\begin{align}\label{P1}
P(s,\chi^{(N)})=&\prod_{P_1\in\Aq}\left(1-\frac{1}{|P_1|_2^s}+O_q\left(|P_1|_2^{-3\sigma/2}\right)\right)\notag\\
=&\zeta_q^{-1}(2s)\prod\limits_{P_1\in\Aq}\left(1+O_q\left(|P_1|_2^{-3\sigma/2}\right)\right),
\end{align}
where $\zeta_q(s)=\frac{1}{1-q^{1-s}}$ is the zeta function for field $\mathbb{F}_q(t)$. 
For an even-degree prime $P_1$, its two factors over $\mathbb F_{q^2}$ have half its degree; this accounts for the exponent $3\sigma/2$ in \eqref{P1}. The error product converges absolutely in the stated region. On each compact subset $K$ of $\Re(s)>1/2$, the finite Euler factors satisfy
\begin{align}\label{P2}
 \prod_{\substack{P_1\mid N\\P_1\in\mathcal A_q}}
 \left|1-|P_1|_2^{-s}\right|^{-1}
 \ll_{K,q,\varepsilon}|N|_2^\varepsilon.
\end{align}
Indeed, factors at sufficiently large primes are bounded by an arbitrarily small power of the prime norm, and the remaining finitely many primes contribute a constant.

Taking $u=q^{-2s}$ and $v=q^{-w}$, we have 
\begin{align}\label{A}
    A_4(u,v)=\sum_{N\in\Aq}v^{\deg N}\frac{L_{q^2}(u, \chi^{(N)})}{L_{q^2}(u^2, \chi^{2(N)})}P(u,\chi^{(N)})\prod_{\substack{P_1|N\\ P_1\in\Aq}}\left(1-u^{\deg P_1}\right)^{-1},
\end{align}
where $$L_{q^2}(u,\chi^{(N)})=\sum\limits_{F\in\Atwo}\chi_F(N)u^{\deg F}$$ 
and 
$$P(u, \chi^{(N)})=\prod\limits_{P_1\in\Aq}(1-u^{\deg P_1}\prod\limits_{_{P_2|P_1}}(1+\chi_{P_2}(N)u^{\deg P_2})^{-1}).$$ It is easy to see that $P(u, \chi^{(N)})$ is analytical for $|u|<q^{-1}$. 

Let $K$ be a compact subset of $|u|<q^{-1}$. Quartic reciprocity identifies $\chi^{(N)}$ with the quartic residue character modulo $N$, so it is principal precisely when $N$ is a fourth power. If $N$ is not a fourth power, passage to the inducing primitive character, Lemma~\ref{LLH}, and the finite Euler factor estimate used in \eqref{P2} give
\begin{align}\label{L}
 \sup_{u\in K}|L_{q^2}(u,\chi^{(N)})|
 \ll_{K,q,\varepsilon}q^{\varepsilon\deg N}.
\end{align}
The same argument includes the quadratic characters arising when $N$ is a square but not a fourth power. The conductor of the inducing character divides $N$; the omitted Euler factors contribute at most $q^{\varepsilon\deg N}$, after relabeling $\varepsilon$.

If $N=M^4$, then
\[
 L_{q^2}(u,\chi^{(N)})
 =\frac{1}{1-q^2u}\prod_{\substack{\pi\mid M\\\pi\in\mathcal A_{q^2}}}
 (1-u^{\deg\pi}).
\]
Consequently $(u-q^{-2})L_{q^2}(u,\chi^{(N)})$ is holomorphic and is also bounded by $O_{K,q,\varepsilon}(q^{\varepsilon\deg N})$. For every $N$, the inverse of $L_{q^2}(u^2,\chi^{2(N)})$ is represented by an absolutely convergent Euler product on $|u|<q^{-1}$ and is uniformly bounded on $K$. The product $P(u,\chi^{(N)})$ is likewise uniformly bounded there, since its local factors are $1+O_K(|u|^{\deg P_1})$.

Applying these bounds to \eqref{A}, with $\varepsilon$ relabeled, yields
\begin{align}\label{S1-regularized-bound}
 |(u-q^{-2})A_4(u,v)|
 \ll_{K,q,\varepsilon}\sum_{N\in\mathcal A_q}
 (q^\varepsilon|v|)^{\deg N},\qquad u\in K.
\end{align}
Choosing $\varepsilon$ small on each compact subset proves normal convergence of the regularized series and hence holomorphy of $(u-q^{-2})A_4(u,v)$ on
\begin{align}
 S_1=\{(u,v):|u|<q^{-1},\ |v|<q^{-1}\}.
\end{align}

  \subsection{Residue at $u=q^{-2}$}\label{PZ}
  We now analyze the possible poles at $u=q^{-2}$. Let $\zeta_{q^2}(u)=\frac{1}{1-q^2u}$ be the zeta function for $\mathbb{F}_{q^2}(T)$ in the variable $u$. The possible pole of $A_4(u,v)$ comes from the term $N$ being a $4$-th power. We have
  \begin{align}     
  &\Res_{u=q^{-2}}\frac{A_4(u,v)}{u^{g/3+2}}\notag\\
  =&\Res_{u=q^{-2}}\frac{1}{u^{g/3+2}}\sum_{N\in\Aq}v^{4\deg N}\frac{\zeta_{q^2}(u)}{\zeta_{q^2}(u^2)}\prod_{\substack{P_2\in\Atwo\\P_2|N}}\frac{1-u^{\deg P_2}}{1-u^{2\deg P_2}}\prod_{P_1 \not| N}(1-u^{\deg P_1}\prod_{\substack{P_2\in\Atwo\\ P_2|P_1}}(1+u^{\deg P_2})^{-1})\notag\\
  =&\Res_{u=q^{-2}}\frac{\zeta_{q^2}(u)}{u^{g/3+2}\zeta_{q^2}(u^2)}P(u)Z(u,v),
  \end{align}
where 
\begin{equation}
    P(u)=\prod_{P_1\in\Aq}(1-u^{\deg P_1}\prod_{\substack{P_2\in\Atwo\\ P_2|P_1}}(1+u^{\deg P_2})^{-1})
\end{equation}
and
\begin{equation}
    Z(u,v)=\sum_{N\in\Aq}v^{4\deg N}\prod_{\substack{P_2\in\Atwo\\P_2|N}}(1+u^{\deg P_2})^{-1}\prod_{\substack{P_1\in\Aq\\ P_1|N}}(1-u^{\deg P_1}\prod_{\substack{P_2\in\Atwo\\ P_2|P_1}}(1+u^{\deg P_2})^{-1})^{-1}.
\end{equation}
  Taking $v=q^{-\frac{1}{2}}$, we then have
  \begin{flalign}
  \Res_{u=q^{-2}}\frac{A_4(u,q^{-\frac{1}{2}})}{u^{g/3+2}}=q^{\frac{2g}{3}}(1-q^{2})P(q^{-2})Z(q^{-2},q^{-\frac{1}{2}}).
  \end{flalign}

We also record that the leading constant is strictly positive. At $u=q^{-2}$, the local factor of $P(u)$ corresponding to a prime $P\in\mathcal A_q$ of degree $d$ is
\[
 p_P=\begin{cases}
 (1+q^{-2d})^{-1},&d\text{ odd},\\[2pt]
 \dfrac{1+2q^{-d}}{(1+q^{-d})^2},&d\text{ even}.
 \end{cases}
\]
These factors are positive and satisfy $p_P=1+O(q^{-2d})$. Since $\sum_Pq^{-2\deg P}<\infty$, their product converges to a strictly positive number. In the defining series for $Z(q^{-2},q^{-1/2})$, the factor associated with each distinct prime divisor $P\mid N$ is
\[
 b_P=\begin{cases}1,&d\text{ odd},\\(1+2q^{-d})^{-1},&d\text{ even}.\end{cases}
\]
Thus $0<b_P\le1$, and
\[
 1\le Z(q^{-2},q^{-1/2})
 =\sum_{N\in\mathcal A_q}q^{-2\deg N}\prod_{P\mid N}b_P
 \le\sum_{N\in\mathcal A_q}q^{-2\deg N}
 =\frac{1}{1-q^{-1}}<\infty.
\]
Consequently,
\begin{equation}\label{positive-main-constant}
 C_q=(q^2-1)P(q^{-2})Z(q^{-2},q^{-1/2})>0.
\end{equation}

 \subsection{The second region}
Write $v=q^{-w}$ and use $L_q(v,\chi)$ for the polynomial in $v$ associated with $L_q(w,\chi)$. For $\deg F\ge1$, Lemma~\ref{LLH} gives
\[
 |L_q(v,\chi_F)|\ll_{q,\varepsilon}|F|_2^\varepsilon
 \qquad (|v|\le q^{-1/2}).
\]
Consequently, the defining series for $\widehat A_4(u,v)$ converges normally on compact subsets satisfying $|u|<q^{-2}$ and $|v|\le q^{-1/2}$.
For $|v|>q^{-1/2}$, Lemma~\ref{fe} gives
\begin{align*}
 \widehat A_4(u,v)
 &=\frac{1}{qv^2}\frac{1-v}{1-1/(qv)}
 \sum_{\substack{F\in\mathcal H_{q^2}\\\deg F\ge1,\ (F,F^\sigma)=1}}
 \omega(\chi_F)L_q(1/(qv),\overline{\chi_F})(quv^2)^{\deg F}.
\end{align*}
Since $|\omega(\chi_F)|=1$, on each compact set $K$ in this range,
\[
 |L_q(v,\chi_F)|\ll_{K,q,\varepsilon}
 (q^{1+\varepsilon}|v|^2)^{\deg F}.
\]
There are at most $q^{2n}$ polynomials $F$ of degree $n$, so the latter series converges normally when $|u|^{1/2}|v|<q^{-3/2}$. Combining the two estimates (and choosing $\varepsilon$ small relative to the strict inequalities) gives holomorphy of $\widehat A_4$ on the open region
\begin{align}
 S_2=\{(u,v): |u|<q^{-2},\ |u|^{1/2}|v|<q^{-3/2}\}.
\end{align}

  \section{The convex hull of all regions}
  
  \subsection{The third convergence region}
We use $\kappa_q$ from Lemma~\ref{quartic-restriction}. The restriction on $F$ means that $F$ has no nonconstant divisor in $\mathbb F_q[T]$, equivalently $(F,F^\sigma)=1$ for square-free $F$.
Combining the functional equation in Lemma~\ref{fe} with \eqref{quartic-gauss-transfer} gives, initially where the dual series converges absolutely,
\begin{align}\label{dual-with-sign}
 A_4(s,w)
 &=q^{2w-1}\frac{1-q^{-w}}{1-q^{w-1}}
 \sum_{\substack{F\in\mathcal H_{q^2}\\(F,F^\sigma)=1}}
 \frac{\kappa_q^{\deg F}G_{q^2}(1,F)
 L_q(1-w,\overline{\chi_F})}{|F|_2^{s+w}}.
\end{align}
For the constant term $F=1$, this identity is interpreted using the zeta functional identity
\[
 \frac1{qv^2}\frac{1-v}{1-1/(qv)}\frac1{1-1/v}
 =\frac1{1-qv},\qquad v=q^{-w};
\]
it does not use the primitive-character Gauss sum formula for $F=1$.

We now apply M\"obius inversion to remove the restriction on $F$, and expand the dual $L$-function for $\Re(w)<0$. Writing the polynomial before inversion as $M=DF$, the sign splits as
\[
 \kappa_q^{\deg M}=\kappa_q^{\deg D}\kappa_q^{\deg F}.
\]
This step uses the degree formula in \eqref{quartic-gauss-transfer}, rather than attempting to factor $\chi_M(M^\sigma)$ for the intermediate polynomials.
For $(D,F)=1$ and $(DF,N)=1$, quartic reciprocity over $\mathbb F_{q^2}$ and Lemma~\ref{Gauss1} yield
\begin{align*}
 \overline{\chi_{DF}}(N)G_{q^2}(1,DF)
 &=\overline{\chi_D}(N)G_{q^2}(1,D)
 \chi_F(D)^2\overline{\chi_F}(N)G_{q^2}(1,F)\\
 &=G_{q^2}(N,D)G_{q^2}(ND^2,F).
\end{align*}
Thus
\begin{align}\label{dual-expanded-sign}
 A_4(s,w)
 &=q^{2w-1}\frac{1-q^{-w}}{1-q^{w-1}}
 \sum_{N\in\mathcal A_q}\frac1{|N|^{1-w}}
 \sum_{\substack{D\in\mathcal A_q\\(D,N)=1}}
 \frac{\mu(D)\kappa_q^{\deg D}G_{q^2}(N,D)}{|D|_2^{s+w}}
 \notag\\[-2pt]
 &\qquad\qquad\times
 \sum_{\substack{F\in\mathcal A_{q^2}\\(F,DN)=1}}
 \frac{\kappa_q^{\deg F}G_{q^2}(ND^2,F)}{|F|_2^{s+w}}.
\end{align}
We have removed the square-free restriction on the last sum because, when $(F,ND)=1$, its Gauss sum vanishes unless $F$ is square-free.

Retain the untwisted definitions
\begin{align}
 H(Q,z)&=\sum_{\substack{F\in\mathcal A_{q^2}\\(F,Q)=1}}
 G_{q^2}(Q,F)z^{\deg F},\\
 C_4(u,v)&=\sum_{N\in\mathcal A_q}(qv)^{-\deg N}
 \sum_{\substack{D\in\mathcal A_q\\(D,N)=1}}
 \mu(D)G_{q^2}(N,D)(uv^2)^{\deg D}H(ND^2,uv^2).
\end{align}
With $u=q^{-2s}$ and $v=q^{-w}$, equation~\eqref{dual-expanded-sign} becomes
\begin{equation}\label{A-C-signed}
 A_4(u,v)=\frac1{qv^2}\frac{1-v}{1-1/(qv)}C_4(\kappa_q u,v).
\end{equation}
In particular,
\begin{align*}
 C_4(\kappa_q u,v)
 &=\sum_{N\in\mathcal A_q}(qv)^{-\deg N}
 \sum_{\substack{D\in\mathcal A_q\\(D,N)=1}}
 \mu(D)G_{q^2}(N,D)(\kappa_q uv^2)^{\deg D}
 H(ND^2,\kappa_q uv^2).
\end{align*}
For $q\equiv7\pmod8$ the argument is unchanged; for $q\equiv3\pmod8$ it is $-u$. The coefficient variable of $A_4$ itself remains $u$.

Replacing $q$ by $q^2$ in Lemma~\ref{upper2}, and putting $\sigma=\Re(s)$ and $\tau=\Re(w)$, gives
\begin{align}\label{H-corrected}
 |H(ND^2,\kappa_q uv^2)|
 &\ll_{K,q,\varepsilon}|ND^2|_2^{\frac12(\frac32-\sigma-\tau)+\varepsilon}
 \notag\\
 &\ll_{K,q,\varepsilon}
 \big(q^{3/2+\varepsilon}|u|^{1/2}|v|\big)^{\deg N+2\deg D},
\end{align}
after relabeling $\varepsilon$. Here $K$ is a compact subset of the range
\[
 q^{-3/2}<|u|^{1/2}|v|\le q^{-1},
\]
avoiding the possible poles $(\kappa_q uv^2)^4=(uv^2)^4=q^{-10}$. Since $|\kappa_q|=1$, this change of argument preserves all modulus conditions and degree exponents. In particular, the endpoint $|u|^{1/2}|v|=q^{-3/2}$ is not supplied by the cited corollary.
Using \eqref{gauss-D-bound}, namely $|G_{q^2}(N,D)|\le q^{\deg D}$ for square-free $D$ coprime to $N$, we obtain, on such compact sets with $|v|>1$,
\begin{align*}
 |C_4(\kappa_q u,v)|
 &\ll_{K,q,\varepsilon}
 \sum_{N\in\Aq}\frac{1}{|qv|^{\deg N}}
 \sum_{\substack{D\in\Aq\\(D,N)=1}}
 q^{\deg D}|uv^2|^{\deg D}
 \big(q^{3/2+\varepsilon}|u|^{1/2}|v|\big)^{\deg N+2\deg D}\\
 &\ll_{K,q,\varepsilon} \sum_{N\in\Aq}
 \big(q^{1/2+\varepsilon}|u|^{1/2}\big)^{\deg N}
 \sum_{D\in\Aq}
 \big(q^{4+2\varepsilon}|u|^2|v|^4\big)^{\deg D}\\
 &=\sum_{n\ge0}\big(q^{3/2+\varepsilon}|u|^{1/2}\big)^n
 \sum_{d\ge0}\big(q^{5+2\varepsilon}|u|^2|v|^4\big)^d.
\end{align*}
Thus the estimates establish normal convergence in the interior
\[
 S_3^\circ=\{|u|<q^{-3},\ |u|^{1/2}|v|<q^{-5/4},\ |v|>1\}
 \cap\{q^{-3/2}<|u|^{1/2}|v|\le q^{-1}\}.
\]
The intended region including the endpoint is
\[
 S_3=\{|u|<q^{-3},\ |u|^{1/2}|v|<q^{-5/4},\ |v|>1\}
 \cap\{q^{-3/2}\le |u|^{1/2}|v|\le q^{-1}\}.
\]
Its upper bound $q^{-1}$ is redundant because $q^{-5/4}<q^{-1}$. In real-part coordinates the three regions are
\begin{align*}
 B_1&=\{\sigma>1/2,\ \tau>1\},\\
 B_2&=\{\sigma>1,\ \sigma+\tau>3/2\},\\
 B_3&=\{\sigma>3/2,\ 5/4<\sigma+\tau\le3/2\}.
\end{align*}
Here $B_j$ denotes the base corresponding to $S_j$ under $u=q^{-2s}$, $v=q^{-w}$. Geometrically, $B=B_1\cup B_2\cup B_3$ is connected and open: indeed,
\[
 B_2\cup B_3=\{\sigma>1,\ \sigma+\tau>3/2\}
 \cup\{\sigma>3/2,\ \sigma+\tau>5/4\}.
\]
The convex hull is
\begin{align}\label{real-hull}
 \widehat B=\{\sigma>1/2,\ \sigma+\tau>5/4,\ 10\sigma+8\tau>13\}.
\end{align}
The last boundary joins $(1/2,1)$ and $(3/2,-1/4)$. Its image in the multiplicative coordinates is
\begin{align}\label{region-S4}
 S_4=\{(u,v): |u|<q^{-1},\ |u|^{1/2}|v|<q^{-5/4},\ |u|^5|v|^8<q^{-13}\}.
\end{align}

\begin{remark}[Pending analytic step]\label{pending-gluing}
The geometry above does not by itself prove analytic continuation. To apply Theorem~\ref{Bochner}, one still has to establish that
\[
 F(s,w)=(q^{-2s}-q^{-2})\widehat A_4(s,w)
\]
is one holomorphic function on $T(B)$, including complex neighborhoods crossing $\Re(s+w)=3/2$ with $\Re(s)>3/2$, and that the local expressions agree. The cited bound only establishes the strict range on one side. The required endpoint-neighborhood argument remains to be supplied. The continuation and contour argument below, and hence the stated asymptotic and its non-vanishing consequence, are conditional on this pending step. No conclusion about a secondary term is drawn here.
\end{remark}

\subsection{Completion of the proof of Theorem~\ref{main}, conditional on the pending analytic step}\label{conditional-proof}
Assume the holomorphy and compatibility specified in Remark~\ref{pending-gluing}. Bochner's theorem then extends $F$ holomorphically to $T(\widehat B)$.
We explain why this extension descends to a single-valued function of $u$ and $v$.
On the original tube, $F$ is invariant under the translations
\[
 (s,w)\mapsto(s+\pi i/\log q,w),\qquad
 (s,w)\mapsto(s,w+2\pi i/\log q),
\]
as follows first from the defining series and then from the identity theorem. Each translate of the holomorphic extension is defined on the same connected tube $T(\widehat B)$ and agrees with $F$ on the original tube. The identity theorem therefore preserves both periods throughout the extension. These periods generate the fibers of
\[
 (s,w)\mapsto(q^{-2s},q^{-w}).
\]
Consequently, the extension descends to a single-valued holomorphic function on $S_4\cap(\mathbb C^*)^2$. Dividing by $u-q^{-2}$ gives a single-valued meromorphic continuation of $\widehat A_4(u,v)$ there, with at most a simple pole at $u=q^{-2}$. At $u=0$ with $v=q^{-1/2}$, use the original convergent power series.

Fix $v=q^{-1/2}$. Region~\eqref{region-S4} permits $|u|<q^{-9/5}$. Choose
\[
 0<\delta<1/5,\qquad 0<r<q^{-2}<R=q^{-9/5-\delta}<q^{-9/5}.
\]
For admissible genera with $g/3+1$ an integer, Cauchy's coefficient formula and the residue theorem, with both circles oriented counterclockwise, give
\begin{align*}
 \sum_{\substack{\chi\ \mathrm{primitive\ quartic}\\\chi^2\ \mathrm{primitive}\\\operatorname{genus}(\chi)=g}}
 L_q(1/2,\chi)
 &=\frac{1}{2\pi i}\oint_{|u|=r}
 \frac{\widehat A_4(u,q^{-1/2})}{u^{g/3+2}}\,du\\
 &=-\Res_{u=q^{-2}}\frac{\widehat A_4(u,q^{-1/2})}{u^{g/3+2}}
 +\frac{1}{2\pi i}\oint_{|u|=R}
 \frac{\widehat A_4(u,q^{-1/2})}{u^{g/3+2}}\,du.
\end{align*}
By \eqref{A-positive}, $\widehat A_4$ and $A_4$ have the same residue at $u=q^{-2}$. In particular, the residue in Section~\ref{PZ} retains the factor $1-q^2$, whereas its contribution to the main term has the factor $q^2-1$.
The continued function is bounded on the fixed compact circle $|u|=R$, so
\begin{align*}
 \left|\frac{1}{2\pi i}\oint_{|u|=R}
 \frac{\widehat A_4(u,q^{-1/2})}{u^{g/3+2}}\,du\right|
 &\ll_{q,\delta}R^{-g/3-1}\\
 &\ll_{q,\delta}q^{(3/5+\delta/3)g}.
\end{align*}
Choosing $\delta$ sufficiently small in terms of $\varepsilon$ yields, subject to Remark~\ref{pending-gluing},
\[
 \sum_{\substack{\chi\ \mathrm{primitive\ quartic}\\\chi^2\ \mathrm{primitive}\\\operatorname{genus}(\chi)=g}}
 L_q(1/2,\chi)
 =q^{2g/3}(q^2-1)P(q^{-2})Z(q^{-2},q^{-1/2})
 +O_{q,\varepsilon}(q^{(3/5+\varepsilon)g}).
\]

 \subsection{Proof of Corollary~\ref{non}}
Let
\begin{align*}
    N_4=\#\{\mbox{$\chi$ quartic of genus g, $\chi$ and $\chi^2$ primitive:}\ L_q(\frac{1}{2},\chi)\not=0\}.
\end{align*}
Theorem~\ref{main}, together with \eqref{positive-main-constant} and a choice of $0<\varepsilon<1/15$ in its error term, gives
\begin{align*}
    q^{\frac{2}{3}g}\ll \sum_{\substack{\chi\ primitive\ quartic\\ \chi^2 primitive\\ genus(\chi)=g}}|L_q(\frac{1}{2}, \chi)|.
\end{align*}
Combining this with \eqref{equation} and Cauchy--Schwarz inequality, we have
\begin{align*}
    q^{\frac{2}{3}g}\ll \sum_{\substack{\chi\ primitive\ quartic\\ \chi^2 primitive\\ genus(\chi)=g}}|L_q(\frac{1}{2}, \chi)|\ll\left(\sum_{\substack{F\in\Htwo \\ (F,F^\sigma)=1}}|L_q(\frac{1}{2}, \chi_{F})|^2\right)^{\frac{1}{2}}N_4^{\frac{1}{2}}.
\end{align*}
Using Perron's formula and the Lindel\"of Hypothesis (Lemma~\ref{LLH}) , we have 
\begin{align}
    \sum_{\substack{F\in\Htwo \\ (F,F^\sigma)=1}}|L_q(\frac{1}{2}, \chi_{F})|^2&\ll q^{4\varepsilon (\frac{g}{3}+1)}\sum_{\substack{F\in\Htwo \\ (F,F^\sigma)=1}} 1\\
    &\ll q^{4\varepsilon (\frac{g}{3}+1)}\frac{1}{2\pi i}\oint_{|u|=r}\sum_{\substack{F\in\mathcal{H}_{q^2} \\ (F,F^\sigma)=1}}u^{\deg F}\frac{du}{u^{g/3+2}}.
\end{align}
Meanwhile, we have
\begin{align}\label{size}
    \sum_{\substack{F\in\mathcal{H}_{q^2} \\ (F,F^\sigma)=1}}u^{\deg F}=&\sum_{\substack{D\in\Aq}}\mu(D)u^{\deg D}\sum_{\substack{(F,D)=1\notag\\ F\in\mathcal{H}_{q^2}}}u^{\deg F}\\
    =&\prod_{P_1\in\Aq}(1-u^{\deg P_1}\prod_{P_2|P_1}(1+u^{\deg P_2})^{-1})\frac{\zeta^{-1}_{q^2}(u^2)}{\zeta_{q^2}^{-1}(u)}.
\end{align}
Then the residue of the above function gives
\begin{align}
    \sum_{\substack{F\in\Htwo \\ (F,F^\sigma)=1}} 1\ll q^{\frac{2g}{3}}
\end{align}
and
\begin{align*}
     \sum_{\substack{F\in\Htwo \\ (F,F^\sigma)=1}}|L_q(\frac{1}{2}, \chi_{F})|^2\ll q^{4\varepsilon (\frac{g}{3}+1)}q^{\frac{2}{3}g}.
\end{align*}
Finally, before relabeling $4\varepsilon$ as $\varepsilon$, we have
\begin{align*}
    N_4\gg q^{\frac{4}{3}g}q^{-4\varepsilon (\frac{g}{3}+1)}q^{-\frac{2}{3}g}=q^{\frac{2}{3}g}q^{-4\varepsilon (\frac{g}{3}+1)}.
\end{align*}

  \section*{Acknowledgement}
The authors would like to express their sincere gratitude to Professor Peng Gao for proposing the subject of this study and for providing invaluable guidance and insightful recommendations throughout the research process.
%\newpage

%\section*{}
%During the preparation of this work the authors used ChatGPT in order to refine the paper. After using this tool, the authors reviewed and edited the content as needed and take full responsibility for the content of the publication.

\bibliographystyle{plain}
\bibliography{ref}

@article{2019The,
  title={The mean values of cubic ${L}$-functions over function fields},
  author={ David, C.  and  Florea, A.  and M. Lalin},
  journal={Algebra \&amp; Number Theory},
  year={2019},
}

@article{ZHENG2017460,
title = {Davenport--{H}asse's theorem for polynomial Gauss sums over finite fields},
journal = {Journal of Number Theory},
volume = {180},
pages = {460-473},
year = {2017},
issn = {0022-314X},
doi = {https://doi.org/10.1016/j.jnt.2017.04.005},
url = {https://www.sciencedirect.com/science/article/pii/S0022314X17301919},
author = {Z. Zheng}
}

@article{BAIER2010879,
title = {Mean values with cubic characters},
journal = {Journal of Number Theory},
volume = {130},
number = {4},
pages = {879-903},
year = {2010},
issn = {0022-314X},
doi = {https://doi.org/10.1016/j.jnt.2009.11.007},
url = {https://www.sciencedirect.com/science/article/pii/S0022314X10000119},
author = {S. Baier and M. P. Young}
}

@article{GAO2024125,
title = {First moment of central values of some primitive Dirichlet ${L}$-functions with fixed order characters},
journal = {Journal of Number Theory},
volume = {261},
pages = {125-142},
year = {2024},
issn = {0022-314X},
doi = {https://doi.org/10.1016/j.jnt.2024.02.007},
url = {https://www.sciencedirect.com/science/article/pii/S0022314X24000593},
author = {P. Gao and L. Zhao}
}

@article{https://doi.org/10.1112/plms/s3-54.2.193,
author = {Patterson, S. J.},
title = {The Distribution of General {G}auss Sums and Similar Arithmetic Functions at Prime Arguments},
journal = {Proceedings of the London Mathematical Society},
volume = {s3-54},
number = {2},
pages = {193-215},
doi = {https://doi.org/10.1112/plms/s3-54.2.193},
url = {https://londmathsoc.onlinelibrary.wiley.com/doi/abs/10.1112/plms/s3-54.2.193},
eprint = {https://londmathsoc.onlinelibrary.wiley.com/doi/pdf/10.1112/plms/s3-54.2.193},
year = {1987}
}

@article{Boc,
      author={Bochner, S.},
       title={A theorem on analytic continuation of functions in several
  variables},
        date={1938},
     journal={Ann. of Math. (2)},
      volume={39},
      number={1},
       pages={14-19},
}

@article{Cech1,
      author={{\v C}ech, M.},
       title={The {R}atios conjecture for real {D}irichlet characters and
  multiple {D}irichlet series},
        date={2024},
     journal={Trans. Amer. Math. Soc.},
      volume={377},
       pages={3487-3528},
}

@article{soundararajan2010second,
  title={The second moment of quadratic twists of modular ${L}$-functions},
  author={Soundararajan, K. and Young, M. P},
  journal={Journal of the European Mathematical Society},
  volume={12},
  number={5},
  pages={1097--1116},
  year={2010}
}

@article{2017On,
  title={On the distribution of the rational points on cyclic covers in the absence of roots of unity},
  author={L. Bary-Soroker  and  P. Meisner },
  journal={Mathematika},
  year={2017},
}

@article{1992Theta,
  title={Theta functions on then-fold metaplectic cover of {SL}(2)?the function field case},
  author={ Hoffstein, J. },
  journal={Inventiones mathematicae},
  volume={107},
  number={1},
  pages={61-86},
  year={1992},
}

@article{andrade2012mean,
  title={The mean value of {$L$}(1/2, $\chi$) in the hyperelliptic ensemble},
  author={Andrade, J. C. and Keating, J. P.},
  journal={Journal of Number Theory},
  volume={132},
  number={12},
  pages={2793--2816},
  year={2012},
  publisher={Elsevier}
}

@article{florea2017improving,
  title={Improving the error term in the mean value of in the hyperelliptic ensemble},
  author={Florea, A. M},
  journal={International Mathematics Research Notices},
  volume={2017},
  number={20},
  pages={6119--6148},
  year={2017},
  publisher={Oxford University Press}
}

@article{lin2024one,
  title={One-level density of zeros of {D}irichlet {$L$}-functions over function fields},
  author={Lin, H.},
  journal={Journal of Number Theory},
  volume={258},
  pages={368--413},
  year={2024},
  publisher={Elsevier}
}

@article{berg2024vanishing,
  title={Vanishing of quartic and sextic twists of {$L$}-functions},
  author={Berg, J. and Ryan, N. C and Young, M. P},
  journal={Research in Number Theory},
  volume={10},
  number={1},
  pages={20},
  year={2024},
  publisher={Springer}
}

@article{gao2021moments,
  title={Moments of central values of quartic {D}irichlet {$L$}-functions},
  author={Gao, P. and Zhao, L.},
  journal={Journal of Number Theory},
  volume={228},
  pages={342--358},
  year={2021},
  publisher={Elsevier}
}

@inproceedings{dunn2025quartic,
  title={Quartic Gauss sums over primes and metaplectic theta functions},
  author={Dunn, A.},
  booktitle={2025 Joint Mathematics Meetings (JMM 2025)},
  organization={AMS}
}

@misc{hong2025meanvaluecubiclfuncitons,
      title={Mean value of cubic {$L$}-funcitons with fixed genus}, 
      author={Z. Hong and Z. Fang and Z. Zheng},
      year={2025},
      eprint={2503.17228},
      archivePrefix={arXiv},
      primaryClass={math.NT},
      url={https://arxiv.org/abs/2503.17228}, 
}

@article{david2025nonvanishing,
  title={Nonvanishing of {$L$}--functions associated to fixed order characters over function fields},
  author={David, C. and Florea, A. and Lalin, M.},
  journal={arXiv preprint arXiv:2506.07815},
  year={2025}
}

@article{D2000RANDOM,
title={RANDOM MATRICES, FROBENIUS EIGENVALUES, AND MONODROMY (American Mathematical Society Colloquium Publications 45) By NICHOLAS M. KATZ and PETER SARNAK: 419 pp., ISBN 0 8218 1017 0 (American Mathematical Society, 1998).},
author={D. R. HeathBrown},
journal={Bulletin of the London Mathematical Society},
volume={32},
number={01},
pages={105-123},
year={2000},
}

@article{Katz1999Zeroes,
title={Zeroes of zeta functions and symmetry.},
author={Katz and Nicholas, M. and Sarnak and Peter},
journal={Bulletin (New Series) of the American Mathematical Society},
year={1999},
}

@misc{DIPP2026,
 author={Diaconu, Adrian and Ion, Bogdan and Pa\c{s}ol, Vicen\c{t}iu and Popa, Alexandru A.},
 title={On the second moment and non-vanishing of central values of {Hecke} {$L$}-functions of {$r$}-th order characters},
 year={2026},
 eprint={2607.27131},
 archivePrefix={arXiv},
 primaryClass={math.NT},
 url={https://arxiv.org/abs/2607.27131}
}

\end{document}